\documentclass{amsart}

\usepackage{epsfig}
\usepackage{amsmath, amsthm, amssymb}

\newtheorem{defn}{Definition}[section]

\newtheorem{cor}[defn]{Corollary}
\newtheorem{thm}[defn]{Theorem}

\newtheorem{lemma}[defn]{Lemma}

\newcommand{\be}{\begin{equation}}
\newcommand{\ee}{\end{equation}}
\newcommand{\bea}{\begin{eqnarray}}
\newcommand{\eea}{\end{eqnarray}}
\newcommand{\beas}{\begin{eqnarray*}}
\newcommand{\eeas}{\end{eqnarray*}}

\newcommand{\R}{\mathbb{R}}

\newcommand{\ve}{\varepsilon}

\newcommand{\noi}{\noindent}

\newcommand{\goto}{\rightarrow}
\newcommand{\hsp}{\hspace{.3in}}
\newcommand{\bp}{\begin{proof}}
\newcommand{\ep}{\end{proof}}

\newcommand{\dstyle}{\displaystyle}

\begin{document}

\title{Mixing Times for the Mean-Field Blume-Capel Model via Aggregate Path Coupling}

\author{Yevgeniy Kovchegov}
\address{Department of Mathematics, Oregon State University, Corvallis, OR  97331}
\email{kovchegy@science.oregonstate.edu}

\author{Peter T. Otto}
\address{Department of Mathematics, Willamette University, Salem, OR 97302}
\email{potto@willamette.edu}

\author{Mathew Titus}
\address{Department of Mathematics, Oregon State University, Corvallis, OR  97331}
\email{titusm@math.oregonstate.edu}
\subjclass[2000]{Primary 60J10; Secondary 60K35}

\date{\today}



\begin{abstract}
In this paper we investigate the relationship between the mixing times of the Glauber dynamics of a statistical mechanical system with its thermodynamic equilibrium structure. For this we consider the mean-field Blume-Capel model, one of  the simplest statistical mechanical models that exhibits the following intricate phase transition structure: within a two dimensional parameter space there exists a curve at which the model undergoes a second-order, continuous phase transition, a curve where the model undergoes a first-order, discontinuous phase transition, and a tricritical point which separates the two curves. We determine the interface between the regions of slow and rapid mixing. In order to completely determine the region of rapid mixing, we employ a novel extension of the path coupling method, successfully proving rapid mixing even in the absence of contraction between neighboring states.
\end{abstract}

\maketitle

\section{Introduction}  

\medskip

The notion of mixing times is fundamental in the study of stochastic processes and their applications. In addition to its purely theoretical value, mixing times are widely used in computational and physical sciences. In computer science, mixing times represent running  times for randomized algorithms such as Metropolis-Hastings that are responsible for providing high accuracy solutions to a large class of problems where classical algorithms have impractical running times. In physics, one important use of mixing times is to determine the feasibility of simulating statistical mechanical systems via so called Glauber dynamics. 

\vskip 0.1 in
\noindent
An important question of mixing times for statistical mechanical models is its relationship with the thermodynamic phase transition structure of the system.  More specifically, as a system undergoes an equilibrium phase transition with respect to some parameter; e.g. temperature, how does the corresponding mixing times behave?  This question has been well studied with recent rigorous results for one class of statistical mechanical models, namely those that undergo second-order, continuous phase transitions like the famous Ising model \cite{LPW} \cite{LLP} \cite{DLP}.  For these models, it has been shown that the mixing times undergo a transition at precisely the thermodynamic phase transition point.  However, for models that exhibit the other type of phase transition: first-order, discontinuous; e.g. Potts model with $q > 2$ \cite{Wu}, \cite{CET} and the Blume-Capel model \cite{Blu}, \cite{BEG}, \cite{Cap1}-\cite{Cap3}, \cite{EOT} with weak interaction, as far as well know, this question has not been previously addressed.  In this paper, we prove that the mixing time transition does not coincide with the thermodynamic equilibrium phase transition for a model that exhibits a first-order, discontinuous phase transition.

\vskip 0.1 in
\noindent
First-order, discontinuous phase transitions are more intricate than their counterparts, which makes rigorous analysis of these models traditionally more difficult.  Furthermore, the more complex phase transition structure causes the models to fall outside the scope of standard mixing time techniques including the so-called `path coupling' method \cite{BD}.  The path coupling argument is a powerful tool used to prove rapid mixing for various Markov chains, including the Glauber dynamics for statistical mechanical models that exhibit a second-order, continuous phase transition. However, due to the intricacy of a first-order, discontinuous phase transition, the standard path coupling argument can not be applied in all cases because contraction between couplings of all neighboring states does not exist.

\vskip 0.1 in
\noindent
In this paper we prove the mixing time rates for the mean-field Blume-Capel (BC) model \cite{Blu}, \cite{BEG}, \cite{Cap1}-\cite{Cap3}, a statistical mechanical spin model ideally suited for the analysis of the relationship between the thermodynamic equilibrium behavior and mixing times due to its intricate phase transition structure.  Specifically, the phase diagram of the BC model includes a curve at which the model undergoes a second-order, continuous phase transition, a curve where the model undergoes a first-order, discontinuous phase transition, and a tricritical point which separates the two curves.  

\vskip 0.1 in
\noindent
As mentioned above, in a subset of the first-order, discontinuous phase transition region of the BC model, the standard path coupling method does not apply as it does throughout the second-order, continuous phase transition region.  Therefore, we apply an innovative new method we call {\it aggregate path coupling} in order to prove rapid mixing for the BC model in the entire first-order, discontinuous phase transition region.  While the standard path coupling method assumes contraction between {\it every} pair of states, our extended approach loosens this condition and requires contraction for only certain pairs of states.  This new method is not restricted to only the BC model and can be applied to other statistical mechanical models that undergo a first-order, discontinuous phase transition and moreover is general enough to be a new tool in the theory of mixing times.

\vskip 0.1 in
\noindent
The paper is organized as follows. The mean-field Blume-Capel model is introduced in Section \ref{bc}. There we describe the equilibrium phase transition structure of the model  and state the large deviation principle of the magnetization.  In Section \ref{mt}, we define the Glauber dynamics for the BC model and the notion of mixing times of Markov chains.  In Section \ref{pc}, we introduce the path coupling for the BC model and prove the form of the mean coupling distance.  In Section \ref{rapid1}, we use the standard path coupling to prove rapid mixing in the single phase region where the BC model undergoes a second-order, continuous phase transition.  In Section \ref{rapid2}, we modify the path coupling argument to prove rapid mixing in the region where the BC model undergoes a first-order, discontinuous phase transition. The slow mixing results are stated and  proved in Section \ref{slow}.

\medskip

\section{Mean-field Blume-Capel model} \label{bc}

\medskip

Statistical mechanical models are defined in terms of the Hamiltonian function.  For the mean-field Blume-Capel model, the Hamiltonian function on the configuration space $\Omega^n = \{ -1, 0, 1\}^n$ is defined by 
\[ H_{n, K}(\omega) = \sum_{j=1}^n \omega_j^2 - \frac{K}{n} \left( \sum_{j=1}^n \omega_j \right)^2 \]
for configurations $\omega = (\omega_1, \ldots, \omega_n)$. Here $K$ represents the interaction strength of the model.  Then for inverse temperature $\beta$, the mean-field Blume-Capel model is defined by the sequence of probability measures
\[ P_{n, \beta, K}(\omega) = \frac{1}{Z_n(\beta,K)} \exp \left[-\beta H_{n,K}(\omega) \right] \]
where $Z_n(\beta, K)= \sum_{\omega \in \Omega^n} \exp [-\beta H_{n, K}(\omega)] $ is the normalizing constant called the {\it partition function}.

\vskip 0.1 in
\noindent
In \cite{EOT}, using large deviations theory \cite{EHT}, the authors proved the phase transition structure of the BC model.  The analysis of the $P_{n, \beta, K}$ was facilitated by expressing it
 in the form of a Curie-Weiss (mean-field Ising)-type model.  This is done by absorbing 
the noninteracting component of the Hamiltonian into the product measure $P_n$ that assigns the probability $3^{-n}$ to each $\omega \in \Omega^n$, obtaining
\be
\label{eqn:rewritecanon}
P_{n,\beta,K}(d\omega) = 
\frac{1}{\tilde{Z}_n(\beta,K)} \cdot \exp\!\left[ n \beta K 
\!\left(\frac{S_n(\omega)}{n} \right)^2 \right] P_{n,\beta}(d\omega)
\ee
In this formula $S_n(\omega)$ equals
the total spin $\sum_{j=1}^n \omega_j$,
$P_{n,\beta}$ is the product measure on $\Omega^n$ with identical one-dimensional marginals
\be
\label{eqn:rhobeta}
\rho_\beta(d \omega_j) = 
\frac{1}{Z(\beta)} \cdot \exp(-\beta \omega_j^2) \, \rho(d \omega_j),
\ee
$Z(\beta)$ is the normalizing constant 
$\int_\Lambda \exp(-\beta \omega_j^2) \rho(d \omega_j) = 1 + 2 e^{-\beta}$,
and $\tilde{Z}_n(\beta,K)$ is the normalizing constant
$[Z(\beta)]^n/Z_n(\beta,K)$. 

\vskip 0.1 in
\noindent
Although $P_{n,\beta,K}$ has the form of a Curie-Weiss (mean-field Ising) model
when rewritten as in (\ref{eqn:rewritecanon}), 
it is much more complicated because of the $\beta$-dependent
product measure $P_{n,\beta}$ and the presence of the parameter $K$. 
These complications introduce new features to the BC model described above that are not present in the Curie-Weiss model \cite{Ellis}.

\vskip 0.1 in
\noindent
The starting point of the analysis of the 
phase-transition structure of the BC model is the large deviation
principle (LDP) satisfied by the spin per site or {\it magnetization} $S_n/n$ with respect to $P_{n,\beta,K}$.  
In order to state the form of the rate function, we introduce the cumulant generating function $c_\beta$ of the measure
$\rho_\beta$ defined in (\ref{eqn:rhobeta}); for $t \in \R$ this function is defined by 
\be
\label{eqn:cbeta}
c_\beta(t) = \log \int_\Lambda \exp (t\omega_1) \, \rho_\beta (d\omega_1) = 
\log \!\left[ \frac{1+e^{-\beta}(e^t+e^{-t})}{1+2e^{-\beta}} \right] \nonumber 
\ee
We also introduce the Legendre-Fenchel transform of $c_\beta$, which is defined for $z \in [-1,1]$ by
\[
J_\beta(z) = \sup_{t \in \R} \{tz - c_\beta(t)\}
\]
and is finite for $z \in [-1,1]$.
$J_\beta$ is the rate function in Cram\'{e}r's theorem, which 
is the LDP for $S_n/n$ with respect to the product measures
$P_{n,\beta}$ \cite[Thm.\ II.4.1]{Ellis} and is one of the components of
the proof of the LDP for $S_n/n$ with respect to the BC model $P_{n,\beta,K}$.
This LDP stated in the next theorem is proved in Theorem 3.3 in \cite{EOT}.

\begin{thm}
\label{thm:ldppnbetak}  
For all $\beta > 0$ and $K > 0$, with respect to $P_{n,\beta,K}$,
$S_n/n$ satisfies the LDP on $[-1,1]$ with exponential speed $n$ and rate function
\[
I_{\beta,K}(z) = J_\beta(z) - \beta K z^2 - \inf_{y \in \R}
\{J_\beta(y) - \beta K y^2 \}.
\]
In other words, for any closed subset $F$,
\be
\label{eqn:upperldp}
\limsup_{n \goto \infty} \frac{1}{n} \log P_{n, \beta, K} \{ S_n/n \in F \} \leq - I_{\beta, K}(F) 
\ee
and for any open subset $G$, 
\be
\label{lowerldp} 
\liminf_{n \goto \infty} \frac{1}{n} \log P_{n, \beta, K} \{ S_n/n \in G \} \geq - I_{\beta, K}(G)
\ee
where $I_{\beta,K}(A) = \inf_{z \in A} I_{\beta,K}(z)$.
\end{thm}
\noi
The LDP in the above theorem implies that those $z \in [-1,1]$
satisfying $I_{\beta,K}(z) > 0$ have an exponentially small probability
of being observed as $n \goto \infty$. Hence we define 
the set of equilibrium macrostates by
\[
\tilde{\mathcal{E}}_{\beta, K} = \{z \in [-1,1] : I_{\beta,K}(z) = 0\}.
\]

\vskip 0.1 in
\noindent
For $z \in \R$ we define
\be
\label{eqn:gbetak} 
G_{\beta,K}(z) = \beta K z^2 - c_\beta(2\beta K z)
\ee
and as in \cite{EMO1} and \cite{EMO2} refer to it as the {\it free energy functional} of the model.  The calculation of the zeroes of $I_{\beta,K}$ --- equivalently, the global minimum points
of $J_{\beta,K}(z) - \beta K z^2$ --- is greatly facilitated by the following observations
made in Proposition 3.4 in \cite{EOT}: 
\begin{enumerate}
\item 
The global minimum points of 
$J_{\beta,K}(z) - \beta K z^2$ coincide with the global minimum points of $G_{\beta,K}$,
which are much easier to calculate. 

\item 
The minimum values $\min_{z \in \R}\{J_{\beta,K}(z) - \beta K z^2\}$ 
and $\min_{z \in \R}\{G_{\beta,K}(z)\}$ coincide.
\end{enumerate}
Item 1 gives the alternate characterization that 
\be
\label{eqn:ebetak}
\tilde{\mathcal{E}}_{\beta, K} = \{z \in [-1,1] : z \mbox{ minimizes } G_{\beta,K}(z)\}.
\ee

\begin{figure}[h]
\begin{center}
\includegraphics[height=1.6in]{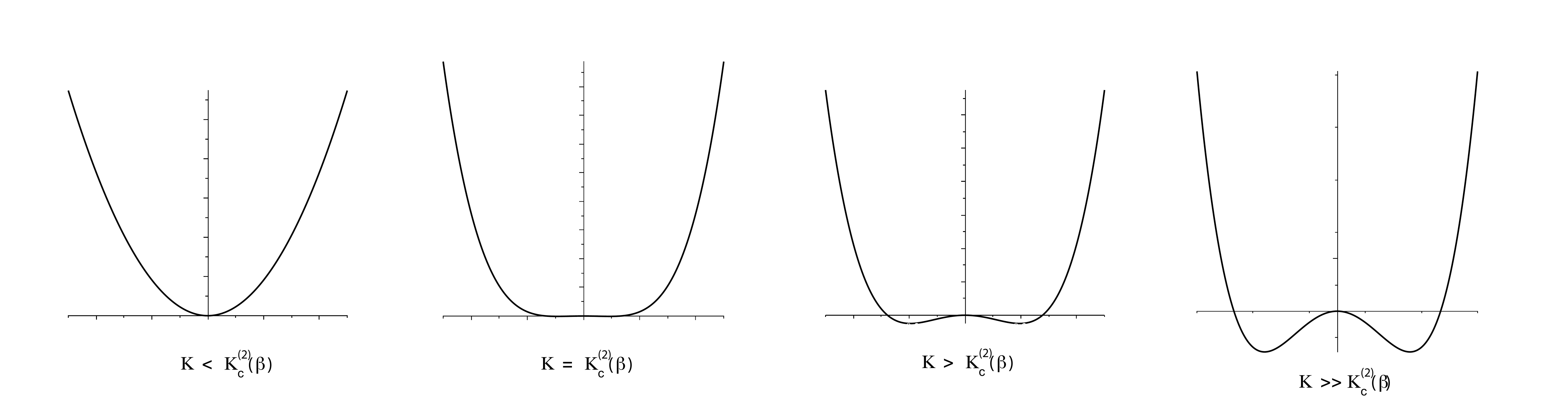}
\caption{\footnotesize The free-energy functional $G_{\beta, K}$ for $\beta \leq \beta_c$} \label{continuous}
\end{center}
\end{figure}

\vskip 0.1 in
\noindent
The free energy functional $G_{\beta, K}$ exhibits two distinct behaviors depending on whether $\beta \leq \beta_c = \log 4$ or $\beta > \beta_c$.  In the first case, the behavior is similar to the Curie-Weiss (mean-field Ising) model.  Specifically, there exists a critical value $K_c^{(2)}(\beta)$ defined in (\ref{eqn:kcbeta}) such that for $K < K_c^{(2)}(\beta)$, $G_{\beta, K}$ has a single minimum point at $z = 0$.  At the critical value $K = K_c^{(2)}(\beta)$, $G_{\beta, K}$ develops symmetric non-zero minimum points and a local maximum point at $z =0$.  This behavior corresponds to a continuous, second-order phase transition and is illustrated in Figure \ref{continuous}.

\begin{figure}[h]
\begin{center}
\includegraphics[height=2.75in]{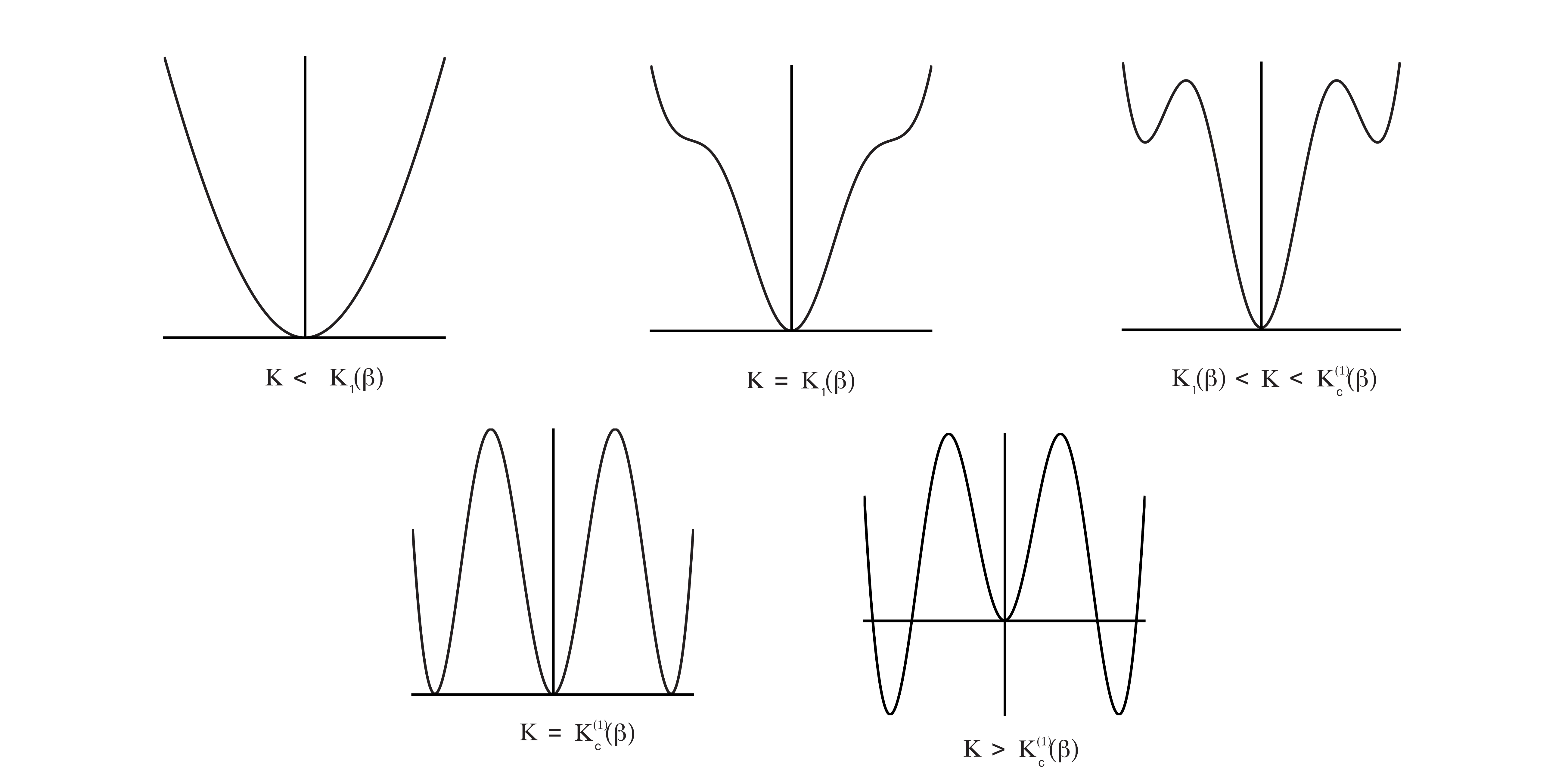}
\caption{\footnotesize The free-energy functional $G_{\beta, K}$ for $\beta > \beta_c$} \label{discont}
\end{center}
\end{figure}

\vskip 0.1 in
\noindent
On the other hand, for $\beta > \beta_c$, $G_{\beta, K}$ undergoes two transitions at the values denoted by $K_1(\beta)$ and $K_c^{(1)}(\beta)$.  For $K < K_1(\beta)$, $G_{\beta, K}$ again possesses a single minimum point at $z = 0$.  At the first critical value $K_1(\beta)$, $G_{\beta, K}$ develops symmetric non-zero local minimum points in addition to the global minimum point at $z=0$.  These local minimum points are referred to as {\it metastable states} and we refer to $K_1(\beta)$ as the {\it metastable critical value}.  This value is defined implicitly in Lemma 3.9 of \cite{EOT} as the unique value of $K$ for which there exists a unique $z >0$ such that
\[  G_{\beta, K_1(\beta)}'(z) = 0 \hsp \mbox{and} \hsp G_{\beta, K_1(\beta)}''(z) = 0 \] 
As $K$ increases from $K_1(\beta)$ to $K_c^{(1)}(\beta)$, the local minimum points decrease until at $K=K_c^{(1)}(\beta)$, the local minimum points reach zero and $G_{\beta, K}$ possesses three global minimum points.  Therefore, for $\beta > \beta_c$, the BC model undergoes a phase transition at $K=K_c^{(1)}(\beta)$, which is defined implicitly in \cite{EOT}.  Lastly, for $K > K_c^{(1)}(\beta)$, the symmetric non-zero minimum points drop below zero and thus $G_{\beta, K}$ has two symmetric non-zero global minimum points.  This behavior corresponds to a discontinuous, first-order phase transition 
 and is illustrated in Figure \ref{discont}.

\vskip 0.1 in
\noindent
In the next two theorems, the 
structure of $\tilde{\mathcal{E}}_{\beta, K}$ corresponding to the behavior of $G_{\beta, K}$ just described is stated which depends on the relationship between $\beta$ and the
critical value $\beta_c = \log 4$.
We first describe $\tilde{\mathcal{E}}_{\beta, K}$ for 
$0 < \beta \leq \beta_c$ and then for $\beta > \beta_c$.  
In the first case $\tilde{\mathcal{E}}_{\beta, K}$ undergoes a continuous
bifurcation as $K$ increases through the critical value $K_c^{(2)}(\beta)$
defined in (\ref{eqn:kcbeta}); physically, this bifurcation
corresponds to a second-order phase transition.  The following theorem 
is proved in Theorem 3.6 in \cite{EOT}.

\begin{thm}
\label{thm:secondorder} 
For $0 < \beta \leq \beta_c$, we define
\be
\label{eqn:kcbeta}
K_c^{(2)}(\beta) = \frac{1}{2\beta c''_\beta(0)} =
\frac{e^\beta + 2}{4\beta}.
\ee
For these values of $\beta$, $\mathcal{E}_{\beta, K}$ has the following structure.

{\em(a)} For $0 < K \leq K_c^{(2)}(\beta)$,
$\tilde{\mathcal{E}}_{\beta,K} = \{0\}$.

{\em(b)} For $K > K_c^{(2)}(\beta)$, there exists 
${z}(\beta,K) > 0$ such that
$\tilde{\mathcal{E}}_{\beta,K} = \{\pm z(\beta,K) \}$.

{\em(c)} ${z}(\beta,K)$ is
a positive, increasing, continuous function for $K > K_c^{(2)}(\beta)$, and
as $K \goto (K_c^{(2)}(\beta))^+$, $z(\beta,K) \goto 0$. 
Therefore, $\tilde{\mathcal{E}}_{\beta,K}$ exhibits a continuous bifurcation
at $K_c^{(2)}(\beta)$.
\end{thm}
\noi
For $\beta \in (0,\beta_c)$, the curve $(\beta,K_c^{(2)}(\beta))$ is the curve of second-order 
critical points.   As we will see in a moment, for $\beta \in (\beta_c,\infty)$
the BC model also has a curve
of first-order critical points, which we denote by $(\beta,K_c^{(1)}(\beta))$.

\vskip 0.1 in
\noindent
We now describe $\tilde{\mathcal{E}}_{\beta, K}$ for 
$\beta > \beta_c$.  In this case $\tilde{\mathcal{E}}_{\beta, K}$ undergoes a discontinuous
bifurcation as $K$ increases through an implicitly defined critical value.
Physically, this bifurcation
corresponds to a first-order phase transition.  The following theorem 
is proved in Theorem 3.8 in \cite{EOT}.  

\begin{figure}[t]
\begin{center}
\includegraphics[height=2.7in]{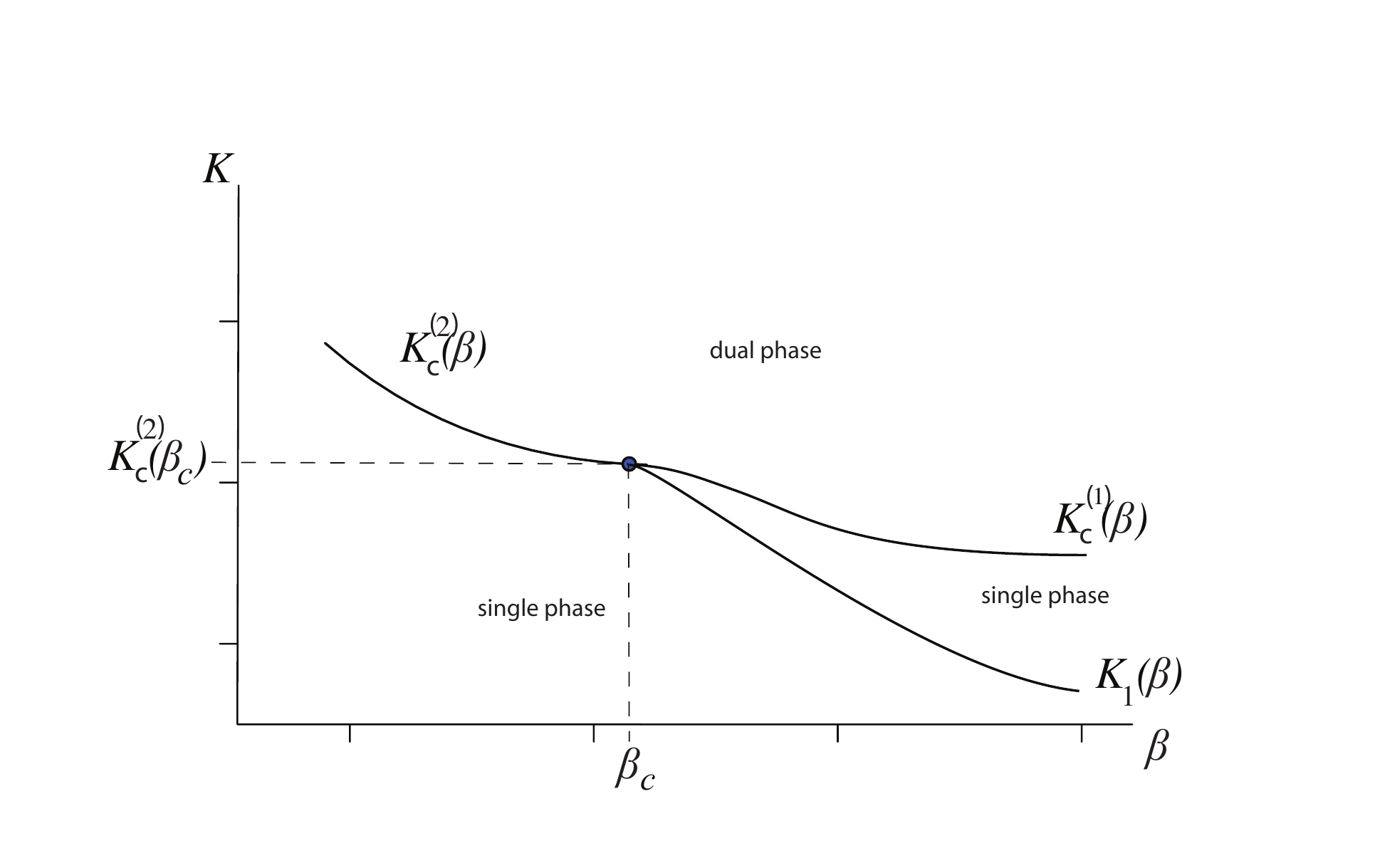}
\caption{\footnotesize Equilibrium phase transition structure of the mean-field Blume-Capel model} \label{phase}
\end{center}
\end{figure}

\begin{thm}
\label{thm:firstorder} 
For all $\beta > \beta_c $, $\tilde{\mathcal{E}}_{\beta, K}$ has
the following structure in terms of the quantity $K_c^{(1)}(\beta)$ defined implicitly for $\beta > \beta_c$ on page {\em 2231} of
{\em \cite{EOT}}.

{\em(a)} For $0 < K < K_c^{(1)}(\beta)$,
$\tilde{\mathcal{E}}_{\beta,K} = \{0\}$.

{\em(b)} There exists $z(\beta,K_c^{(1)}(\beta)) > 0$
such that $\tilde{\mathcal{E}}_{\beta,K_c^{(1)}(\beta)} =
\{0,\pm z(\beta,K_c^{(1)}(\beta))\}$.

{\em(c)} For $K > K_c^{(1)}(\beta)$ 
there exists $z(\beta,K) > 0$
such that $\tilde{\mathcal{E}}_{\beta,K} =
\{\pm z(\beta,K)\}$.

{\em(d)} $z(\beta,K)$
is a positive, increasing, continuous function for $K \geq K_c^{(1)}(\beta)$, and 
as $K \goto K_c^{(1)}(\beta)^+$, $z(\beta,K) \goto 
z(\beta,K_c^{(1)}(\beta)) > 0$.  Therefore,
$\tilde{\mathcal{E}}_{\beta,K}$ exhibits a discontinuous bifurcation
at $K_c^{(1)}(\beta)$.
\end{thm}
\noi
The phase diagram of the BC model is depicted in Figure \ref{phase}.  The LDP stated in Theorem \ref{thm:ldppnbetak} implies the following weak convergence result used in the proof of rapid mixing in the first-order, discontinuous phase transition region.  It is part (a) of Theorem 6.5 in \cite{EOT}.

\begin{thm}
\label{thm:weakconv}
For $\beta$ and $K$ for which $\tilde{\mathcal{E}}_{\beta, K} = \{ 0 \}$, 
\[ P_{n, \beta, K} \{ S_n/n \in dx \} \Longrightarrow \delta_0 \hsp \mbox{as} \ \ n \goto \infty. \]
\end{thm}

\vskip 0.1 in
\noindent
We end this section with a final result that was not included in the original paper \cite{EOT} but will be used in the proof of the slow mixing result for the BC model.  The result states that not only do the global minimum point of $G_{\beta, K}$ and $I_{\beta, K}$ coincide, but so do the local minimum points.

\begin{lemma}
\label{lemma:equalmin}
A point $\tilde{z}$ is a local minimum point of $G_{\beta, K}$ if and only if it is a local minimum point of $I_{\beta, K}$.
\end{lemma}
\bp
Assume that $\tilde{z}$ is a local minimum point of $G_{\beta, K}$.  Then $\tilde{z}$ is a critical point of $G_{\beta, K}$ which implies that $\tilde{z} = c_\beta '(2 \beta K \tilde{z})$.  By the theory of Legendre-Fenchel transforms, $J_\beta '(z) = (c_\beta ')^{-1}(z)$ and thus 
\[ I_{\beta, K}'(\tilde{z}) = J_{\beta} '(\tilde{z}) - 2 \beta K \tilde{z} = (c_\beta ')^{-1}(\tilde{z}) - 2 \beta K \tilde{z} = 0. \]

\vskip 0.1 in
\noindent
Next, since $\tilde{z}$ is a local minimum point of $G_{\beta, K}$, 
\[ G_{\beta, K} ''(\tilde{z}) > 0 \hsp \Longleftrightarrow \hsp c_{\beta}''(2 \beta K \tilde{z}) < \frac{1}{2 \beta K} \]
Therefore,
\[ I_{\beta, K}''(\tilde{z}) = J_{\beta}''(\tilde{z}) - 2 \beta K = \frac{1}{c_{\beta}''(2 \beta K \tilde{z})} - 2 \beta K > 0 \]
and we conclude that $\tilde{z}$ is a local minimum point of $I_{\beta, K}$.  The other direction is obtained by reversing the argument.
\ep

\medskip

\section{Glauber Dynamics and Mixing Times}  \label{mt}

\medskip

The Markov chain for statistical mechanical models studied in this paper is called the Glauber dynamics, one of the most commonly used physical dynamics for these models. See \cite{Bre} for more on Glauber dynamics.  For the BC model, the Glauber dynamics evolve by selecting a vertex $i$ at random and updating the spin at $i$ according to the distribution $P_{n, \beta, K}$, conditioned to agree with the spins at all vertices not equal to $i$.  If the current configuration is $\omega$ and vertex $i$ is selected, then the chance of the spin at $i$ is updated to $+1$ is equal to 
\be
\label{eqn:plustran}
p_{+1} (\omega, i) := \frac{e^{2\beta K \tilde{S}(\omega, i)/n }}{e^{2\beta K \tilde{S}(\omega, i)/n} + e^{\beta-(\beta K)/n} + e^{-2\beta K \tilde{S}(\omega, i)/n} }  \ee
where $\tilde{S}(\omega, i) = \sum_{j\neq i} \omega_j$ is the total spin of the neighboring vertices of $i$.
Similarly, the probabilities of $i$ updating to $0$ and $-1$ are
\be 
\label{eqn:zerotran}
p_{0}(\omega, i) = \frac{e^{\beta-(\beta K)/n}}{e^{2\beta K \tilde{S}(\omega, i)/n} + e^{\beta-(\beta K)/n} + e^{-2\beta K \tilde{S}(\omega, i)/n} } \ee
and
\be 
\label{eqn:minustran}
p_{-1}(\omega, i) = \frac{e^{-2\beta K \tilde{S}(\omega, i)/n} }{e^{2\beta K \tilde{S}(\omega, i)/n} + e^{\beta-(\beta K)/n} + e^{-2\beta K \tilde{S}(\omega, i)/n} }
\ee
$p_{+1} (\omega, i)$ is increasing with respect to $\tilde{S}(\omega, i)$, $p_{-1} (\omega, i)$ is decreasing with respect to $\tilde{S}(\omega, i)$, and $p_{0} (\omega, i)$ is decreasing for $\tilde{S}(\omega, i) > 0$ and increasing for $\tilde{S}(\omega, i) < 0$.

\vskip 0.1 in
\noindent
The mixing time is a measure of the convergence rate of a Markov chain to its stationary distribution and is defined in terms of 
the {\it total variation distance} between two distributions $\mu$ and $\nu$ defined by
\[ \| \mu - \nu \|_{TV} = \sup_{A \subset \Omega} | \mu(A) - \nu(A)| = \frac{1}{2} \sum_{x \in \Omega} | \mu(x) - \nu(x)| \]
Given the convergence of the Markov chain, we define the maximal distance to stationary to be 
\[ d(t) = \max_{x \in \Omega} \| P^t(x, \cdot) - \pi \|_{TV} \]
where $P^t(x, \cdot)$ is the transition probability of the Markov chain starting in configuration $x$ and $\pi$ is its stationary distribution.
Then, given $\ve>0$, the {\it mixing time} of the Markov chain is defined by
\[ t_{mix}(\ve) = \min\{ t : d(t) \leq \ve \} \smallskip \]
In general, rigorous analysis of mixing times is very difficult and the proof of exact mixing time asymptotics of even some basic chains remain elusive. See \cite{LPW} for a detailed survey on the theory of mixing times. 

\vskip 0.2 in
\noindent
Rates of mixing times are generally categorized into two groups: {\it rapid mixing} which implies the mixing time exhibits polynomial growth with respect to the system size, and {\it slow mixing} which implies that the mixing time grows exponentially with the system size.  Determining where a model undergoes rapid mixing is of major importance, as it is in this region that the application of the dynamics is physically feasible.  

\vskip 0.2 in
\noindent
A classical tool in proving rapid mixing for Markov chains defined on graphs, including the Glauber dynamics of statistical mechanical models, is the path coupling technique \cite{BD}.  It will be shown that this technique can be directly applied for the BC model in the second-order, continuous phase transition region but fails in a subset of the first-order, discontinuous phase transition region.  For the latter region, we develop an extension of the path coupling method we call {\it aggregate path coupling} to prove rapid mixing.  The path coupling method for the BC model is introduced in the next section.

\medskip

\section{Path Coupling}  \label{pc}

\medskip

We begin by setting up the path coupling for the Glauber dynamics of the mean-field Blume-Capel model.  Define the path metric $\rho$ on $\Omega^n = \{-1, 0, 1\}^n$ by 
\be
\label{eqn:pathmetric} 
\rho(\sigma, \tau) = \sum_{j=1}^n \mathbf{1}\{ \sigma_j \neq \tau_j \}, 
\ee
the number of sites at which the configurations $\sigma$ and $\tau$ differ.  

\vskip 0.2 in
\noindent
Let $\sigma$ and $\tau$ be two configurations with $\rho(\sigma, \tau) = 1$; i.e. $\sigma$ and $\tau$ are neighboring configurations.  The spins of $\sigma$ and $\tau$ agree everywhere except at a single vertex $i$. Assume that $\sigma_i < \tau_i$.  We next describe the path coupling $(X, Y)$ of one step of the Glauber dynamics starting in configuration $\sigma$ with one starting in configuration $\tau$.
Pick a vertex $k$ uniformly at random.  We use a single random variable as the common source of noise to update both chains, so the two chains agree as often as possible.  In particular, let $U$ be a uniform random variable on $[0,1]$ and set
\[ X(k) = \left\{ \begin{array}{rl} -1 & \mbox{if \ $0 \leq U \leq p_{-1}(\sigma, k)$} \\ 0 & \mbox{if \ $p_{-1}(\sigma, k) < U \leq p_{-1}(\sigma, k) + p_{0}(\sigma, k)$} \\ +1 & \mbox{if \ $p_{-1}(\sigma, k) + p_{0}(\sigma, k) < U \leq 1$} \end{array} \right. \]
and
\[ Y(k) = \left\{ \begin{array}{rl} -1 & \mbox{if \ $0 \leq U \leq p_{-1}(\tau, k)$} \\ 0 & \mbox{if \ $p_{-1}(\tau, k) < U \leq p_{-1}(\tau, k) + p_{0}(\tau, k)$} \\ +1 & \mbox{if \ $p_{-1}(\tau, k) + p_{0}(\tau, k) < U \leq 1$} \end{array} \right. \]
Set $X(j) = \sigma_j$ and $Y(j) = \tau_j$ for $j \neq k$.

\vskip 0.2 in
\noindent
Since $\sigma_i < \tau_i$, for all $j \neq i$, $\tilde{S}(\sigma, j) < \tilde{S}(\tau, j)$ and thus
\[ p_{+1}(\tau, k) > p_{+1}(\sigma, k) \ \ \ \mbox{and} \ \ \ p_{-1}(\tau, k) < p_{-1}(\sigma, k) \]
The path metric $\rho$ on the coupling above takes on the following possible values.
\[ \rho (X,Y) = \left\{ \begin{array}{lll} 0 & \mbox{if} & k = i \\ 1 & \mbox{if} & \mbox{$k \neq i$ and both chains updates the same}
\\ 2 & \mbox{if} & \mbox{$k \neq i$ and the chains update differently} \end{array} \right.  \]

\noi
The application of the path coupling technique to prove rapid mixing is dependent on whether the mean coupling distance with respect to the path metric $\rho$, denoted by $\mathbb{E}_{\sigma, \tau} [ \rho(X,Y)] $, contracts over all pairs of neighboring configurations.

\vskip 0.2 in
\noindent
In the lemma below and following corollary, we derive a working form for the mean coupling distance.

\begin{lemma}
\label{lemma:meanpath}
Let $\rho$ be the path metric defined in {\em (\ref{eqn:pathmetric})} and $(X,Y)$ be the path coupling where $X$ and $Y$ start in neighboring configurations $\sigma$ and $\tau$.  Define 
\be
\label{eqn:varphi} 
\varphi_{\beta, K}(x) = \frac{2\sinh (\frac{2\beta K}{n}x)}{ 2\cosh (\frac{2\beta K}{n} x) + e^{\beta - \frac{\beta K}{n}}} 
\ee
Then 
\[ \mathbb{E}_{\sigma, \tau} [ \rho(X,Y)]  = \frac{n-1}{n} + \frac{(n-1)}{n} [\varphi_{\beta, K}(S_n(\tau)) - \varphi_{\beta, K}(S_n(\sigma))] + O\left(\frac{1}{n^2} \right) \]
\end{lemma}
\bp
Let $n_{-1}, n_0$ and $n_{+1}$ denote the number of $-1, 0$ and $+1$ spins, respectively, not including the spin at vertex $i$, where the configurations differ, in configuration $\sigma$.  Note that $n_{-1} + n_0 + n_{+1} = n-1$.  

\vskip 0.2 in
\noindent
Define $\varepsilon(-1)$ to be the probability that $X$ and $Y$ update differently when the chosen vertex $k \neq i$ is a $-1$ spin.  Similarly, define $\varepsilon(0)$ and $\varepsilon(+1)$.  Then
\beas
\mathbb{E}_{\sigma, \tau} [ \rho(X,Y)] & = & \frac{n_{-1}}{n} (1 - \varepsilon(-1)) +  \frac{n_{0}}{n} (1 - \varepsilon(0)) +  \frac{n_{+1}}{n} (1 - \varepsilon(+1)) \\
 & & + 2 \left[ \frac{n_{-1}}{n} \varepsilon(-1) +  \frac{n_{0}}{n} \varepsilon(0) +  \frac{n_{+1}}{n} \varepsilon(+1) \right] \\
 & = & \frac{n-1}{n} + \frac{n_{-1}}{n} \varepsilon(-1) +  \frac{n_{0}}{n} \varepsilon(0) +  \frac{n_{+1}}{n} \varepsilon(+1)
\eeas
The probability that $X$ and $Y$ update differently when the chosen vertex $k \neq i$ is a $-1$ spin is given by
\beas 
\ve(-1) & = & \left[ p_{-1}(\sigma, k) - p_{-1}(\tau, k) \right] + \left[(p_{-1}(\sigma, k) + p_{0}(\sigma, k)) - (p_{-1}(\tau, k) + p_{0}(\tau, k)) \right] \\
& = & [ p_{+1}(\tau, k) - p_{+1}(\sigma, k) ] + \left[ p_{-1}(\sigma, k) - p_{-1}(\tau, k) \right] \\
& = & [ p_{+1}(\tau, k) - p_{-1}(\tau, k) ] + \left[ p_{-1}(\sigma, k) -  p_{+1}(\sigma, k) \right] \\
& = & \frac{2\sinh (\frac{2\beta K}{n}(S_n(\tau)+1))}{ 2\cosh (\frac{2\beta K}{n} (S_n(\tau)+1)) + e^{\beta - \frac{\beta K}{n}}} - \frac{2\sinh (\frac{2\beta K}{n}(S_n(\sigma)+1))}{ 2\cosh (\frac{2\beta K}{n} (S_n(\sigma)+1)) + e^{\beta - \frac{\beta K}{n}}} \\
& = & \varphi_{\beta, K} ((S_n(\tau) +1)) - \varphi_{\beta, K} ((S_n(\sigma) +1))\\
& = & \varphi_{\beta, K} (S_n(\tau) ) - \varphi_{\beta, K} (S_n(\sigma))+ O\left(\frac{1}{n^2} \right)
\eeas
Similarly, we have 
$$\ve(0) =\varphi_{\beta, K} (S_n(\tau) ) - \varphi_{\beta, K} (S_n(\sigma))$$
 and 
 $$\ve(+1) = \varphi_{\beta, K} ((S_n(\tau) -1)) - \varphi_{\beta, K} ((S_n(\sigma) -1)) =  \varphi_{\beta, K} (S_n(\tau) ) - \varphi_{\beta, K} (S_n(\sigma))+ O\left(\frac{1}{n^2} \right)$$  
 Thus the proof is complete.
\ep

\noindent
For $c_\beta$ defined in (\ref{eqn:cbeta}),  
\[ \varphi_{\beta,K}(x)= c_\beta ' \left( \frac{2 \beta K}{n} x \right) (1+O(1/n)) \]
Thus we have the following corollary.
\medskip
\begin{cor} \label{cor:meanpath}
Let $\rho$ be the path metric defined in {\em (\ref{eqn:pathmetric})} and $(X,Y)$ be the path coupling where $X$ and $Y$ start in neighboring configurations $\sigma$ and $\tau$.  Then 
\[ \mathbb{E}_{\sigma, \tau} [ \rho(X,Y)]  = \frac{n-1}{n} + \frac{(n-1)}{n} \left[c_\beta ' \left( 2\beta K \frac{S_n(\tau)}{n} \right) - c_\beta ' \left( 2 \beta K \frac{S_n(\sigma)}{n} \right) \right] + O\left(\frac{1}{n^2} \right) \]
\end{cor}

\begin{figure}[t]
\begin{center}
\includegraphics[height=2.5in]{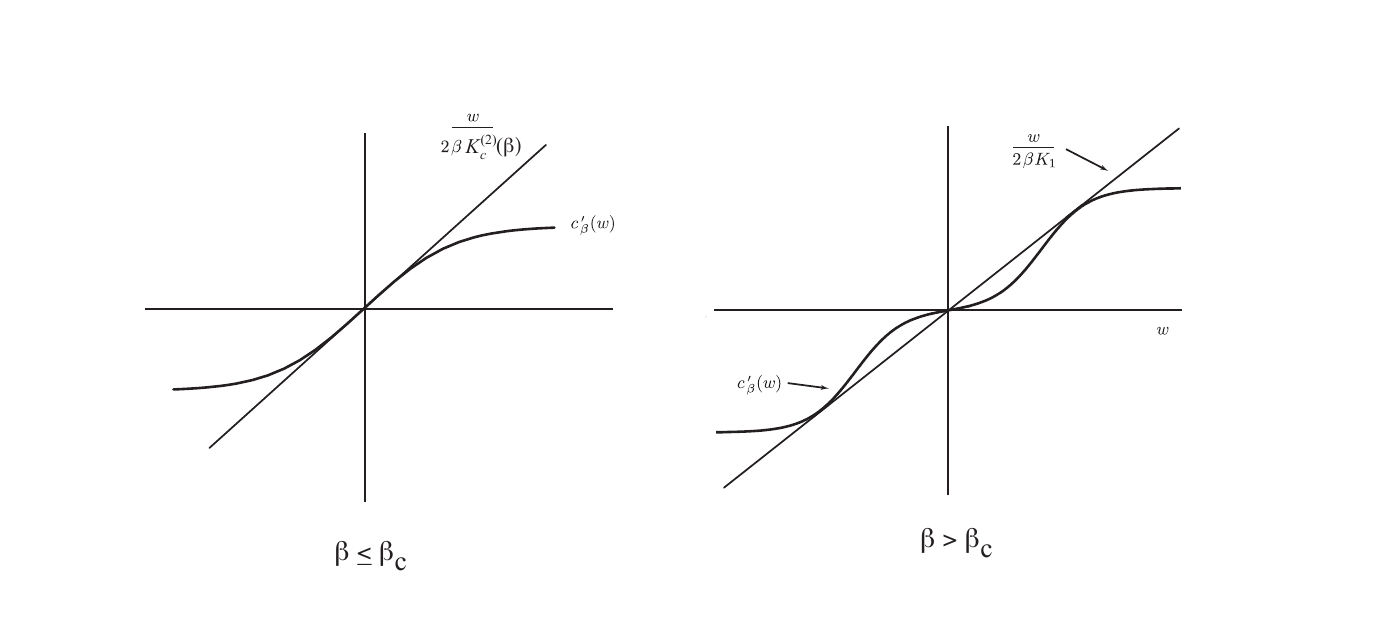}
\caption{\footnotesize Behavior of $c_\beta'(w)$ for large and small $\beta$.} \label{cprime}
\end{center}
\end{figure}

\vskip 0.2 in
\noi
By the above corollary, we conclude that the mean coupling distance of a coupling starting in neighboring configurations contracts; i.e. $\mathbb{E}_{\sigma, \tau}[\rho(X,Y)] < \rho(\sigma, \tau) = 1$, if 
\[ \left[c_\beta ' \left( 2\beta K \frac{S_n(\tau)}{n} \right) - c_\beta ' \left( 2 \beta K \frac{S_n(\sigma)}{n} \right) \right] \approx 2 \beta K \left[ \frac{S_n(\tau)}{n} - \frac{S_n(\sigma)}{n} \right] c_\beta ''\left(2 \beta K \frac{S_n(\sigma)}{n} \right) < \frac{1}{n-1} \]
Since $\sigma$ and $\tau$ are neighboring configurations and $S_n (\tau) > S_n(\sigma)$, this is equivalent to 
\be
\label{eqn:cpp} 
c_\beta ''\left(2 \beta K \frac{S_n(\sigma)}{n} \right) < \frac{1}{2 \beta K} 
\ee
Therefore, contraction of the mean coupling distance, and thus rapid mixing, depends on the concavity behavior of the function $c_\beta'$.  This is also precisely what determines the type of thermodynamic equilibrium phase transition (continuous, second-order versus discontinuous, first-order) that is exhibited by the mean-field Blume-Capel model.  We state the concavity behavior of $c_\beta'$ in the next theorem which is proved in Theorem 3.5 in \cite{EOT}.  The results of the theorem are depicted in Figure \ref{cprime}

\begin{thm} 
\label{thm:concavity}
For $\beta > \beta_c = \log 4$ define 
\be
\label{eqn:wcrit}
w_c(\beta) = \cosh^{-1} \! \left( \frac{1}{2} e^\beta -
4e^{-\beta}  \right) \geq 0.  
\ee 
The following conclusions hold.

{\em (a)} For $0 < \beta \leq \beta_c$, $c_\beta'(w)$ is strictly concave  for
$w > 0$.

{\em (b)} For $\beta > \beta_c$,
$c_\beta'(w)$ is strictly convex for $0 < w < w_c(\beta)$ \ 
and $c_\beta'(w)$ is strictly concave for $w > w_c(\beta)$.
\end{thm} 

\noi
By part (a) of the above theorem, for $\beta \leq \beta_c$, $~~c_\beta''(x)\leq c_\beta''(0) = 1/(2 \beta K_c^{(2)}(\beta))$.  Therefore, by (\ref{eqn:cpp}), the mean coupling distance contracts between {\bf all pairs} of neighboring states whenever $K<K_c^{(2)}(\beta)$. 

\vskip 0.2 in
\noi
By contrast, for $\beta>\beta_c$, we will show that rapid mixing occurs whenever $K < K_1(\beta)$ where $K_1(\beta)$ is the metastable critical value introduced in Section \ref{bc} and depicted in Figure \ref{discont}.   However, since the supremum  $~~\sup_{[-1,1]} c_\beta''(x)>{1 \over 2\beta K_1(\beta)}$, having $K<K_1(\beta)$ is not sufficient for (\ref{eqn:cpp}) to hold. That is, $K<K_1(\beta)$ does not imply the contraction of the mean coupling distance between {\bf all pairs} of neighboring states.
Regardless, we will prove rapid mixing for all $K<K_1(\beta)$ in Section \ref{rapid2} using an extension to the path coupling method that we refer to as {\it aggregate path coupling}.

\vskip 0.2 in
\noi
We now prove the mixing times for the mean-field Blume-Capel model, which varies depending on the parameter values $(\beta, K)$ and their position with respect to the thermodynamic phase transition curves.

\medskip

\section{Rapid Mixing for Continuous Phase Transition Region} \label{rapid1}  

\medskip

We begin by stating the standard path coupling argument used to prove rapid mixing for the mean-field Blume-Capel model in the continuous, second-order phase transition region.  The result is proved in Corollary 14.7 of \cite{LPW}.

\begin{thm}
\label{thm:path}
Suppose the state space $\Omega$ of a Markov chain is the vertex set of a graph with path metric $\rho$.  Suppose that for each edge $\{\sigma,\tau\}$ there exists a coupling $(X,Y)$ of the distributions $P(\sigma, \cdot)$ and $P(\tau, \cdot)$ such that 
\[ \mathbb{E}_{\sigma, \tau}[\rho (X,Y)] \leq \rho(\sigma, \tau) e^{-\alpha} \quad \text{ for some }\alpha>0 \]
Then
\[ t_{mix}(\ve) \leq \left\lceil \frac{-\log (\ve) + \log (\mbox{{\em diam}}(\Omega))}{\alpha} \right\rceil \]
\end{thm}

\vskip 0.2 in
\noi
In this section, we assume $\beta \leq \beta_c$ which implies that the BC model undergoes a continuous, second-order phase transition at $K = K_c^{(2)}(\beta)$ defined in (\ref{eqn:kcbeta}).  By Theorem \ref{thm:concavity}, for $\beta \leq \beta_c$, $c_\beta'(x)$ is concave for $x > 0$.  See the first graph of Figure \ref{cprime} as reference. 
We next state and prove the rapid mixing result for the mean-field Blume-Capel model in the second-order, continuous phase transition regime.

\begin{thm}
\label{thm:mfBC}
Let $t_{mix}(\ve)$ be the mixing time for the Glauber dynamics of the mean-field Blume-Capel model on $n$ vertices and $K_c^{(2)}(\beta)$ be the continuous phase transition curve defined in {\em (\ref{eqn:kcbeta})}.  Then for $\beta \leq \beta_c = \log 4$ and $K < K_c^{(2)}(\beta)$,
\[ t_{mix}(\ve) \leq {n \over \alpha} (\log n + \log (1/\ve)) \]
for any $\alpha \in \left(0, {K_c^{(2)}-K  \over K_c^{(2)}}\right)$ and $n$ large enough.
\end{thm}

\begin{proof}
Let $(X,Y)$ be a coupling of the Glauber dynamics of the BC model that begin in neighboring configurations $\sigma$ and $\tau$ with respect to the path metric $\rho$ defined in (\ref{eqn:pathmetric}).  By Corollary \ref{cor:meanpath} of Lemma \ref{lemma:meanpath}, 
\[
\mathbb{E}_{\sigma, \tau} [ \rho(X,Y)]  = 1 - \left( \frac{1}{n} - \frac{(n-1)}{n} \left[c_\beta ' \left( 2\beta K \frac{S_n(\tau)}{n} \right) - c_\beta ' \left( 2 \beta K \frac{S_n(\sigma)}{n} \right) \right] \right) + O\left(\frac{1}{n^2} \right) \]
Observe that $c_\beta ''$ is an even function and that for $\beta \leq \beta_c$, ${\dstyle \sup_x c_\beta ''(x) = c_\beta ''(0)}$.  Therefore, by the mean value theorem and Theorem \ref{thm:secondorder},
\beas 
\mathbb{E}_{\sigma, \tau} [ \rho(X,Y)] & \leq & 1 - \frac{[1 - (n-1)(2\beta K /n) c_\beta ''(0)]}{n} + O\left(\frac{1}{n^2} \right) \\
& \leq & \exp \left\{- \frac{1-2\beta K c_\beta ''(0)}{n} + O\left( \frac{1}{n^2} \right) \right\} \\
& = &   \exp \left\{\frac{1}{n} \left({K  \over K_c^{(2)}(\beta)}-1 \right) + O \left(\frac{1}{n^2} \right) \right\} \\
& < & e^{-\alpha/n}
\eeas
for any $\alpha \in \left(0, {K_c^{(2)}-K  \over K_c^{(2)}}\right)$ and $n$ sufficiently large.
Thus, for $K < K_c^{(2)}(\beta)$, we can apply Theorem \ref{thm:path} where the diameter of the configuration space of the BC model $\Omega^n$ is $n$ to complete the proof.
\ep

\medskip

\section{Rapid Mixing for Discontinuous Phase Transition Region} \label{rapid2}  

\medskip

Here we consider the region $\beta>\beta_c$, where the mean-field Blume-Capel model undergoes a first-order discontinuous phase transition.  In this region, the function $c_\beta'(x)$ which determines whether the mean coupling distance contracts (Corollary \ref{cor:meanpath}) is no longer strictly concave for $x > 0$ (Theorem \ref{thm:concavity}).  See the second graph in Figure \ref{cprime} for reference.  We will show that rapid mixing occurs whenever $K < K_1(\beta)$ where $K_1(\beta)$ is the metastable critical value defined in Section \ref{bc} and depicted in Figure \ref{discont}.

\vskip 0.2 in
\noi
As shown in Section \ref{pc}, in order to apply the standard path coupling technique of Theorem \ref{thm:path}, we need the inequality (\ref{eqn:cpp}) to hold for all values of $S_n(\sigma)$ and thus $~~\sup_{[-1,1]} c_\beta''(x)<{1 \over 2\beta K}$.  However since $~~\sup_{[-1,1]} c_\beta''(x)>{1 \over 2\beta K_1(\beta)}$, the condition $K<K_1(\beta)$ is not sufficient for the contraction of the mean coupling distance between {\bf all pairs} of neighboring states which is required to prove rapid mixing using the standard path coupling technique stated in Theorem \ref{thm:path}. 

\vskip 0.2 in
\noi 
In order to prove rapid mixing in the region where $\beta>\beta_c$ and $K<K_1(\beta)$, we take advantage of the result in Theorem \ref{thm:weakconv} which states the weak convergence of the magnetization $S_n/n$ to a point-mass at the origin.   Thus, in the coupling of the dynamics, the magnetization of the process that starts at equilibrium will stay mainly near the origin.  As a result, for two starting configurations $\sigma$ and $\tau$, one of which has near-zero magnetization ($S_n(\sigma)/n \approx 0$), the mean coupling distance of a coupling starting in these configurations will be the aggregate of the mean coupling distances between neighboring states along a minimal path connecting the two configurations. Although not all pairs of neighbors in the path will contract, we show that in the {\it aggregate}, contraction between the two configurations still holds.

\vskip 0.2 in
\noi 
In the next lemma we prove contraction of the mean coupling distance in the aggregate and then the rapid mixing result for the mean-field Blume-Capel model is proved in the theorem following the lemma by applying the new aggregate path coupling method.

\begin{lemma}
\label{lemma:aggregate}
Let $(X,Y)$ be a coupling of the Glauber dynamics of the BC model that begin in configurations $\sigma$ and $\tau$, not necessarily neighbors with respect to the path metric $\rho$ defined in {\em (\ref{eqn:pathmetric})}.  Suppose $\beta > \beta_c$ and $K < K_1(\beta)$. Then for any $\alpha \in \left(0,{K_1(\beta)-K \over K_1(\beta)}\right)$ there exists an $\varepsilon>0$ such that, asymptotically as $n \goto \infty$,
\[ \mathbb{E}_{\tau,\sigma} [ \rho(X,Y)] \leq e^{-\alpha/n} \rho(\sigma, \tau) \]
whenever $|S_n(\sigma) |<\varepsilon n$.
\end{lemma}
\bp
Observe that for $\beta > \beta_c$ and $K < K_1(\beta)$, $$|c_\beta '(x)| \leq {|x| \over  2\beta K_1(\beta)} \hsp \mbox{for all} \ x $$ 
We will show that for a given $\alpha' \in \left({1 \over 2\beta K_1(\beta)},{1-\alpha \over 2\beta K} \right)$, there exists $\varepsilon>0$ such that
\begin{equation}\label{cb}
c_\beta '(x)-c_\beta '(x_0) \leq \alpha' (x-x_0)\quad \text{ whenever}~ |x_0|<\varepsilon
\end{equation}
as $c_\beta '(x)$ is a continuously differentiable increasing odd function and $c_\beta '(0)=0$.

\vskip 0.2 in
\noi 
In order to show (\ref{cb}), observe that $c''_{\beta}(0)={1 \over 2\beta K_c^{(2)}(\beta)} < {1 \over 2 \beta K_1(\beta)}$, and since $c''_{\beta}$ is continuous, there exists a $\delta>0$ such that
$$c_\beta ''(x)<\alpha'  \quad \text{ whenever}~ |x|<\delta$$
The mean value theorem implies that
$$c_\beta '(x)-c_\beta '(x_0) <\alpha' (x-x_0)\qquad \forall x_0,x \in (-\delta,\delta)$$
Now, let $\varepsilon={\alpha'-1/(2\beta K_1(\beta)) \over \alpha'+1/(2\beta K_1(\beta))} \delta < \delta$.  Then for any $|x_0|<\varepsilon$ and $|x|\geq \delta$,
$$|c_\beta '(x)-c_\beta '(x_0)| \leq {|x|+|x_0| \over 2\beta K_1(\beta)}
\leq {(1+{\varepsilon/\delta})|x| \over 2\beta K_1(\beta)} 
= {|x-x_0| \over 2\beta K_1(\beta)}\cdot {1+{\varepsilon/\delta} \over |1-x_0/x|} 
\leq {|x-x_0| \over 2\beta K_1(\beta)}\cdot {1+{\varepsilon/\delta} \over 1-\varepsilon/\delta}
=\alpha' |x-x_0|$$

\noindent
Without loss of generality suppose that $S_n(\sigma)< S_n(\tau)$.
Let $(\sigma = x_0, x_1, \ldots, x_r = \tau)$ be a path connecting $\sigma$ to $\tau$ and monotone increasing  in $\rho$ such that $(x_{i-1}, x_i)$ are neighboring configurations.  Here $r=\rho(\sigma,\tau)$. Then by Corollary \ref{cor:meanpath} of Lemma \ref{lemma:meanpath} and (\ref{cb}), we have for $|S_n(\sigma) |<\varepsilon n$ and asymptotically as $n \rightarrow \infty$, 
\beas
 \mathbb{E}_{\tau,\sigma} [ \rho(X,Y)]  & \leq & \sum_{i=1}^r \mathbb{E}_{x_{i-1}, x_i} [ \rho(X_{i-1},X_i)] \\
& = &  \frac{(n-1)}{n} \rho(\sigma, \tau) +  \frac{(n-1)}{n} \left[c_\beta'\left( \frac{2\beta K}{n} S_n(\tau) \right) - c_\beta' \left(\frac{2\beta K}{n} S_n(\sigma) \right) \right]  +\rho(\sigma, \tau) \cdot O\left({1 \over n^2}\right)\\
& \leq & \frac{(n-1)}{n} \rho(\sigma, \tau) +  \frac{(n-1)}{n} (S_n(\tau)-S_n(\sigma)) \frac{2\beta K \alpha'}{n}  +\rho(\sigma, \tau) \cdot O\left({1 \over n^2}\right)\\
& \leq & \rho (\sigma, \tau) \left[ 1 - \left(\frac{1-2\beta K \alpha'}{n} \right)  + O\left({1 \over n^2}\right) \right] \\
& \leq & e^{-\alpha/n} \rho(\sigma, \tau)
\eeas
This completes the proof.
\ep

\begin{thm}
\label{thm:mfBC}
Let $t_{mix}(\ve)$ be the mixing time for the Glauber dynamics of the mean-field Blume-Capel model on $n$ vertices and $K_1(\beta)$ be the metastable critical point.  Then, for $\beta > \beta_c$ and $K < K_1(\beta)$,
\[ t_{mix}(\ve) \leq {n \over \alpha} (\log n + \log (1/\ve)) \]
for any $\alpha \in \left(0,{K_1(\beta)-K \over K_1(\beta)}\right)$ and $n$ large enough.
\end{thm}
\begin{proof}  Let $(X_t, Y_t)$ be a coupling such that $Y_0 \overset{dist}{=} P_{n, \beta, K}$, the stationary distribution.  
For a given $\alpha \in \left(0,{K_1(\beta)-K \over K_1(\beta)}\right)$, let  $\varepsilon$ be as in Lemma \ref{lemma:aggregate}.  Then for sufficiently large $n$,
\beas
\| P^t(X_0, \cdot) - P_{n, \beta, K} \|_{TV} & \leq & P \{ X_t \neq Y_t \}  \\
& = & P \{ \rho(X_t, Y_t ) \geq 1 \} \\
& \leq & \mathbb{E} [ \rho(X_t,Y_t)] \\
& \leq & \mathbb{E} [ \rho(X_t,Y_t)~|~|S_n(Y_{t-1})|<\varepsilon n] \cdot P\{|S_n(Y_{t-1})|<\varepsilon n\}+2nP\{|S_n(Y_{t-1})| \geq \varepsilon n\} \\
& \leq & e^{-\alpha/n} \mathbb{E} [ \rho(X_{t-1},Y_{t-1}); |S_n(Y_{t-1})|<\varepsilon n] \cdot P\{|S_n(Y_{t-1})|<\varepsilon n\}+2nP\{|S_n(Y_{t-1})| \geq \varepsilon n\} \\
& \leq & e^{-\alpha/n} \mathbb{E} [ \rho(X_{t-1},Y_{t-1})] \cdot P\{|S_n(Y_{t-1})|<\varepsilon n\}+2nP\{|S_n(Y_{t-1})| \geq \varepsilon n\} \\
& \vdots & \qquad \vdots\\
& \leq & e^{-\alpha t/n} \mathbb{E} [ \rho(X_0,Y_0)] \cdot \prod_{s=0}^{t-1} P\{|S_n(Y_s)|<\varepsilon n\}+2n \sum_{s=0}^{t-1} P\{|S_n(Y_s)| \geq \varepsilon n\} \\
& \leq & e^{-\alpha t/n} \mathbb{E} [ \rho(X_0,Y_0)] \cdot \left(P_{n, \beta, K}\{|S_n/n|<\varepsilon\}\right)^t+2ntP_{n, \beta, K}\{|S_n/n| \geq \varepsilon\} \\
& \leq & 2n e^{-\alpha t/n} \cdot \left(P_{n, \beta, K}\{|S_n/n|<\varepsilon\}\right)^t+2ntP_{n, \beta, K}\{|S_n/n| \geq \varepsilon\} 
\eeas

\noi
We recall the result in Theorem \ref{thm:weakconv} that for $\beta > \beta_c$ and $K < K_1(\beta)$ 
\[ P_{n, \beta, K}\{ S_n/n \in dx \} \Longrightarrow \delta_0 \hsp \mbox{as $n \goto \infty$.} \]
Moreover, for any $\gamma>1$ and $n$ sufficiently large, the LDP stated in Theorem \ref{thm:ldppnbetak} implies that
\beas
\|P^t(X_0, \cdot) - P_{n, \beta, K} \|_{TV} & \leq & 2n e^{-\alpha t/n} \cdot \left(P_{n, \beta, K}\{|S_n/n|<\varepsilon\}\right)^t+2ntP_{n, \beta, K}\{|S_n/n| \geq \varepsilon\}  \\
&  <  & n e^{-\alpha t/n} \cdot \left(1-e^{-\gamma nI_{\beta,K}(\varepsilon)} \right)^t+tne^{-{n \over \gamma}I_{\beta,K}(\varepsilon)} 
\eeas
This completes the proof. 
\ep

\medskip

\section{Slow Mixing}  \label{slow}

\medskip

We complete our mixing time analysis of the mean-field Blume-Caple model by determining the slow mixing region of the $(\beta, K)$ parameter space for both the continuous, second-order and discontinuous, first-order phase transition regions.  To prove slow mixing for the BC model, we use the bottleneck or Cheeger constant argument which we state in the theorem below.  The theorem is an application of Theorem 7.3 in \cite{LPW} to the BC model.

\begin{thm}
\label{thm:bneck}
For two configurations $\omega$ and $\tau$, define the edge measure $Q$ as follows:
\[ Q(\omega, \tau) = P_{n, \beta, K}(\omega) P(\omega, \tau) \hsp \mbox{and} \hsp Q(A, B) = \sum_{\omega \in A, \tau \in B} Q(\omega, \tau) \] 
Here $P(\omega, \tau)$ is the transition probability corresponding to the Glauber dynamics of the mean-field Blume-Capel model.  The bottleneck ratio of the set $S$ is defined by 
\[ \Phi(S) = \frac{Q(S, S^c)}{P_{n, \beta, K}(S)} \hsp \mbox{and} \hsp \Phi_\ast = \min_{S:  P_{n, \beta, K}(S) \leq \frac{1}{2}} \Phi(S) \]
Under these assumptions, we have
\[ t_{mix} = t_{mix}(1/4) \geq \frac{1}{4 \Phi_\ast} \]
\end{thm}

\medskip

\noi
The next lemma states that bottlenecks, and thus slow mixing, occurs for the BC model whenever the function $G_{\beta, K}$ has a positive minimum point.

\begin{lemma}
\label{lemma:slowmix1}
Suppose the free energy functional $G_{\beta, K}$ has a positive minimum (either local or global) point.  Then there exists a positive constant $b$ and a strictly positive function $r(\beta,K)$ such that
\[ t_{mix} \geq be^{r(\beta,K)n} \]
\end{lemma}
\bp
Suppose $G_{\beta, K}$ has a minimum (either local or global) point at $\tilde{z} > 0$.  Furthermore, let $z'$ be where $G_{\beta, K}$ has a local maximum point such that $0 \leq z' < \tilde{z}$.  Define the bottleneck set
\[ A = \{ \omega : z' < S_n(\omega)/n \leq 1 \} \]
Since $z' \geq 0$, by the symmetry of $G_{\beta, K}$, $P_{n, \beta, K} (A) \leq \frac{1}{2}$.  With respect to the Glauber dynamics, in order to leave the set $A$ in one step of the chain, the current configuration must lie in one of the two boundary sets
\[ A_1 = \left\{ \omega : \frac{S_n(\omega)}{n} = z' + \frac{1}{n} \right\} \hsp \mbox{or} \hsp A_2 = \left\{ \omega : \frac{S_n(\omega)}{n} = z' + \frac{2}{n} \right\} \]
We take $n$ sufficiently large so that $z' + 2/n < \tilde{z}$.

\vskip 0.2 in
\noi
For any $\gamma>1$, the large deviations upper bound (\ref{eqn:upperldp}) implies that the edge measure $Q$ for the bottleneck set $A$ satisfies
\[ Q(A, A^c) \leq P_{n, \beta, K}(A_1) + P_{n, \beta, K}(A_2) < e^{-{n \over \gamma} I_{\beta, K}(z' + \frac{2}{n})} \]
for $n$ large enough.
Moreover, since the LDP implies that $P_{n, \beta, K}(A) > e^{-n \gamma I_{\beta, K}(\tilde{z})}$, the bottleneck ratio satisfies
\[ \Phi(A) < 2 e^{-n [\gamma^{-1}I_{\beta, K}(z' + \frac{2}{n}) - \gamma I_{\beta, K}(\tilde{z})]}, \]
all for a given $\gamma>1$ and $n$ large enough.
Lastly, by Lemma \ref{lemma:equalmin} $\tilde{z}$ is a local minimum point of $I_{\beta, K}$ since it is a local minimum point of $G_{\beta, K}$, and thus $I_{\beta, K}(z') > I_{\beta, K}(\tilde{z})$ and the proof is complete by Theorem \ref{thm:bneck}.
\ep

\noi
Next we state the slow mixing result for the mean-field Blume-Capel model.  The result follows from Lemma \ref{lemma:slowmix1} and the region in the $(\beta, K)$ parameter space at which the free-energy functional $G_{\beta, K}$ possesses a positive minimum point.  These regions were determined in \cite{EOT} and depicted in Figures \ref{continuous} and \ref{discont}.  

\begin{cor}
\label{cor:slowmix}
Let $t_{mix} = t_{mix}(1/4)$ be the mixing time for the Glauber dynamics of the mean-field Blume-Capel model on $n$ vertices.  For {\em (a)} $\beta \leq \beta_c$ and $K>K_c^{(2)}(\beta)$, and {\em (b)} $\beta > \beta_c$ and $K > K_1(\beta)$, there exists a positive constant $b$ and a strictly positive function $r(\beta,K)$ such that
\[ t_{mix} \geq be^{r(\beta,K)n} \]
\end{cor}

\noi
We summarize the mixing time results for the mean-field Blume-Capel model and its relationship to the model's thermodynamic phase transition structure in Figure \ref{mixingphase}.  As shown in the figure, in the second-order, continuous phase transition region ($\beta \leq \beta_c$) for the BC model, the mixing time transition coincides with the equilibrium phase transition.  This is consistent with other models that exhibit this type of phase transition.  On the other hand, this is not the case in the first-order, discontinuous phase transition region ($\beta > \beta_c$) where the mixing time transition occurs below the equilibrium phase transition.

\begin{figure}[t]
\begin{center}
\includegraphics[height=3in]{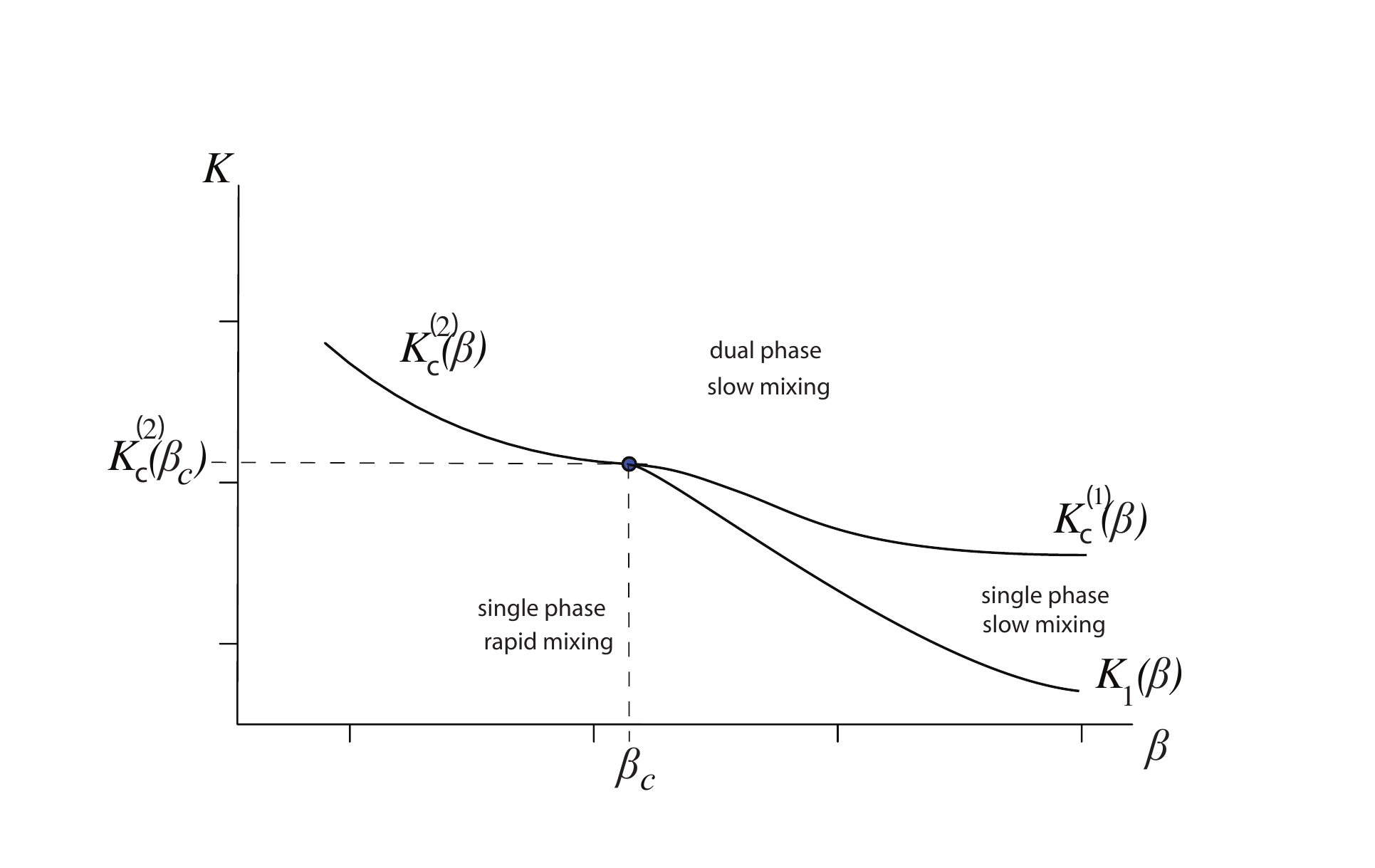}
\caption{\footnotesize Mixing times and equilibrium phase transition structure of the mean-field Blume-Capel model} \label{mixingphase}
\end{center}
\end{figure} 

\bigskip
\bigskip


\bibliographystyle{amsplain}

\end{document}